# Learning Exact Topology of a Loopy Power Grid from Ambient Dynamics

Saurav Talukdar
University of Minnesota
Minneapolis, USA

Deepjyoti Deka
Los Alamos National Lab
Los Alamos, USA

Blake Lundstrom
University of Minnesota
Minneapolis, USA

Michael Chertkov
Los Alamos National Lab
Los Alamos,USA

Murti V. Salapaka
University of Minnesota
Minneapolis, USA

## ABSTRACT

Estimation of the operational topology of the power grid is necessary for optimal market settlement and reliable dynamic operation of the grid. This paper presents a novel framework for topology estimation for general power grids (loopy or radial) using time-series measurements of nodal voltage phase angles that arise from the swing dynamics. Our learning framework utilizes multivariate Wiener filtering to unravel the interaction between fluctuations in voltage angles at different nodes and identifies operational edges by considering the phase response of the elements of the multivariate Wiener filter. The performance of our learning framework is demonstrated through simulations on standard IEEE test cases.

## 1. INTRODUCTION

The power grid comprises a network of transmission lines that enable the flow of electricity between generators and load buses. The grid structure is represented as a graph with nodes denoting buses and edges representing connecting lines. The dynamics of the grid is influenced by the presence of rotating masses that include generator turbines and industrial loads. Moreover, loads that are voltage/frequency dependent with restricted operating regimes affect grid operation. Monitoring the dynamic operation and assessing the small-signal stability of the grid includes efficient estimation of its state variables and topology (set of operational lines). In particular, topology estimation can enable the detection of line failures. However, real-time topology estimation is hindered by the limited presence of real-time line-based measurements including breaker status and flows. The absence of such measurements is more severe in the low and medium voltage lines in the distribution grid [15], where, a majority of the household solar panels are located. Furthermore, the system operator may not also have access to topology of grid areas that are outside its jurisdiction and need to estimate the grid topology indirectly. On the other hand, new devices like phasor measurement units (PMUs) [25], micro-PMUs [30], FNETs [33] that record high-fidelity real-time measurements of nodal states are increasingly being deployed in grid buses. Similarly, smart devices such as air conditioners and electric vehicles often have the ability to monitor nodal voltages for control goals. In this article, we study the problem of estimating the topology of the grid using such voltage measurements collected from grid nodes. Note that due to their high fidelity with sub-second sampling frequency, the measurement samples collected are not independent but represent time-series of nodal voltage dynamics. Our learning approach is thus based on the swing equations [19], which govern the dynamics of nodal voltage angles in the grid. Furthermore, our learning framework does **not** require the values of generator/load parameters, line impedance (susceptance) of permissible lines.

### 1.1 Prior Work

Topology Learning in power grids is a growing area of research. A majority of the prior work have focused on using statistics from static power flow models to learn the topology of radial grids. Approaches include using the signs of inverse covariance matrix of voltages [3], graphical models [5], signature comparison based tests [1], maximum likelihood based tests for line measurements [26] as well as greedy algorithm using voltage second moments [9, 7, 8, 6]. Extending such models for loopy grids is not straightforward as mentioned in [3, 10]. Approximate schemes with good performance are discussed in [21]. However the above mentioned work rely on independent measurements samples that need sufficient time separation to prevent the grid dynamics from introducing correlations.

Here, we relax both the radial assumption on grid topology and the i.i.d. assumption on sample collected. We consider the nodal measurements to arise from the swing dynamics [19] in a loopy power grid. Network topology reconstruction for linear dynamical systems often include active intervention approaches that actively modify inputs (see [32, 23, 27]). However, these approaches are not suitable, as in power grids it is impractical to actively manipulate nodal quantities solely for estimation purposes. In contrast, our approach relies only on passively recorded data. It is shown in [22, 29] that for linear dynamical systems like linearized swing equations, multivariate Wiener filtering based network reconstruction recovers the moral graph [20] of the underlying topology. Here, the reconstructed topology includes spurious edges along with the true edges of the original graph. Other approaches for topology identification using passive observations include using the inverse power spectral density matrix [4] and causally conditioned mutual information [12]. However these are limited to 'directed' linear dynamical systems while power grid dynamics are essentially bidirectional.

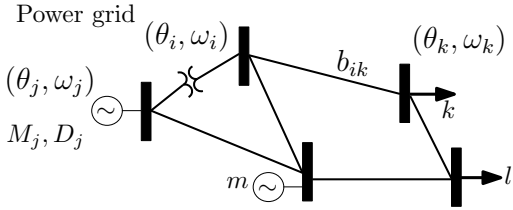

**Figure 1:** Power Grid example: the graph represents lines between the grid nodes/buses. Nodes $j$ and $m$ have generators.

## 1.2 Contribution

Our novel learning framework builds the topology of the grid using time-series measurements of nodal voltage angles in two stages. We consider the linearized swing equations as a model that describes the grid dynamics [19] and learn a multivariate **Wiener filter** [17] for the available nodal phase angle data in the first stage. The sparsity of the multivariate Wiener filter results in a graph that includes the true operational lines/edges of the grid topology as well as additional 'spurious edges' between two-hop neighbors in the grid graph. In the next stage, we develop a **Pruning Step** that separates the true edges from the spurious ones. The pruning step is based on a new result that shows that the Wiener filters associated with spurious edges have constant phase at all frequencies unlike those for true edges. Thus, an exact recovery can be made for all grid graphs, both loopy and radial. To the best of our knowledge, this represents the first work to provide guaranteed topology learning in loopy power grids using data from grid dynamics.

In the next section we present a mathematical model of the grid and the associated swing dynamics under ambient noise. In Section 3, we develop the construction of the Wiener filter based on the swing equations. Following that, in Section 4 we prove the key result of this paper and present our pruning algorithm that recovers the exact topology of the grid. We demonstrate the effectiveness of our algorithm through simulations on standard IEEE test systems in Section 5. Finally conclusions are presented in Section 6.

## 2. MATHEMATICAL FORMULATION

In this section, we provide a mathematical model for structure and dynamics of the power grid.

### 2.1 Swing Equations based model of the Grid

We represent the power grid as a connected undirected graph $\mathcal{G} = \{\mathcal{V}, \mathcal{E}\}$, where $\mathcal{V} = \{1, 2, \cdots, N\}$ is the set of $N$ buses/nodes and $\mathcal{E} = \{(i,j)\}$ is the set of undirected lines/edges, where $(i,j)$ is to be considered as an unordered tuple. Let $b_{ij} > 0$ denote the susceptance of line $(i,j)$ in the grid. The complex voltage at each node $j$, $V_j$ has magnitude and phase angle denoted by $|V_j|$ and $\theta_j$ respectively. The frequency at node $j$ is denoted by $\omega_j$, where $\omega_j = \frac{d\theta_j}{dt}$. The frequency at all nodes is regulated at a constant value $\omega^0 = 60$ Hz (in U.S.A.). For small ambient disturbances in the grid, the temporal dynamics at each node of the grid is represented by the following linearized Swing Equation [18],

$$M_j \ddot{\theta}_j + D_j \dot{\theta}_j = p_j^{(m)} - \sum_{i:(i,j)\in\mathcal{E}} b_{ij}(\theta_j - \theta_i) + p_j, j \in \mathcal{V}, \quad (1)$$

where, $M_j$ denotes the inertia of the rotating mass, $D_j$ denotes the damping, $p_j^{(m)}$ represents the real power injection and $p_j$ is the external disturbance at the node $j$. Here, $p_{ji} = b_{ij}(\theta_j - \theta_i)$ gives the real line flow from node $j$ to $i$. Under equilibrium conditions $\omega_j = \omega^0, \dot{\omega}_j = 0$ and power balance is satisfied at each node $(p_j^{(m)} - \sum_{i:(i,j)\in\mathcal{E}} p_{ji})$. As all terms involving frequency and phase angle are linear, we take the equilibrium point as reference and express dynamics in terms of deviations from the reference. Abusing notation, we use $\omega_j, \theta_j, p_{ji}$ to denote the deviation from their nominal values respectively instead of their true values.

As mentioned earlier, the time-series measurements collected from smart meters comprise of discrete-time samples of phase angles at the grid nodes. Writing Eq. (1) in discrete-time (indexed by $n$) using first order difference for time-derivative of variables we have

$$M_j \frac{\dot{\theta}_j(n+1) - \dot{\theta}_j(n)}{t_s} + D_j \frac{\theta_j(n+1) - \theta_j(n)}{t_s}$$
$$= \sum_{i:(i,j)\in\mathcal{E}} b_{ij}\theta_i(n) - B_j\theta_j(n) + p_j(n) \quad (2)$$

where $B_j = \sum_{i:(i,j)\in\mathcal{E}} b_{ij}$. Subsequently we use its $z$ domain representation, relating the output phase angle at each node in terms of the phase angle of its neighbors and exogenous input $E_j(z)$ as follows:

$$\Theta_j(z) = \sum_{i:(i,j)\in\mathcal{E}} H_{ji}(z)\Theta_i(z) + E_j(z), \forall j \in \mathcal{V}, \quad (3)$$

where, the **transfer function** $H_{ji}(z)$ from $\Theta_i(z)$ to $\Theta_j(z)$ and the input $E_j(z)$ are given by,

$$H_{ji}(z) = \frac{b_{ji}}{S_j(z)}, E_j(z) = \frac{1}{S_j(z)} P_j(z), \quad (4)$$

and, $$S_j(z) = \frac{M_j}{t_s^2}(z-1)^2 + \frac{D_j}{t_s}(z-1) + B_j. \quad (5)$$

The frequency response of $H_{ji}(z)$ is given by $H_{ji}(e^{\hat{j}\omega})$, the magnitude and phase angle of which is denoted as $|H_{ji}(e^{\hat{j}\omega})|$ and $\angle(H_{ji}(e^{\hat{j}\omega}))$ respectively. Next, we discuss the characteristics of the disturbance $p_j$ at each node.

### 2.2 Stochastic Ambient Disturbance

A scalar random process $x(n)$ is *wide sense stationary* (WSS) if its mean $\mu(n) := \mathbb{E}[x(n)]$ is constant and correlation $R_x(s,t) = \mathbb{E}[x(s)x(t)]$ is a function of $s-t$. A vector random processes $x(n) = [x_1(n)\ x_2(n)...x_N(n)]^T$ is said to be *wide sense stationary* (WSS) if $\forall i,j$, $x_i(n), x_j(n)$ are WSS and the cross correlation function $R_{x_i,x_j}(s,t) = \mathbb{E}[x_i(s)y_j(t)]$ is a function of $s-t$. For $x_i(n)$ and $x_j(n)$, the **Cross Power Spectral Density** is given as $\Phi_{x_i,x_j}(z) =$

$\sum_{n=-\infty}^{\infty} R_{x_i,x_j}(n)z^{-n}$. The power spectral density for vector $x(n)$ is given by the matrix $\Phi_x(z)$, where $[\Phi_x(z)](i,j) = \Phi_{x_i x_j}(z)$. If $R_{x_i,x_j}(s,t)$ (subsequently $\Phi_{x_i,x_j}(z)$) is uniformly zero almost everywhere $\forall i, j$, then $x(n)$ is said to be uncorrelated. Note that *the correlation and power spectral density matrices for uncorrelated vector random processes are diagonal. Further the diagonal of the power spectral density matrix is real, even and positive at all frequencies* [24].

**Disturbance Model:** We model the vector of ambient disturbances $p$ at grid nodes by uncorrelated WSS process with zero mean, that is, disturbance at the same node is time-correlated, while at two distinct nodes are uncorrelated. The use of uncorrelated zero mean WSS processes to model ambient disturbances is prevalent in power grids [13, 31]. In the next section, we establish methods to estimate the topology of the grid interconnections.

## 3. WIENER FILTERING BASED RECONSTRUCTION

Consider the swing equations for a grid $\mathcal{G} = \{\mathcal{V}, \mathcal{E}\}$ discussed in Eqs. (3),(4). The output at node $j$, $\theta_j(n)$ is dependent on the state $\theta_i(n)$ of node $i$ through the transfer function $H_{ji}(z)$ and the exogenous input $e_j$ with $z$ transform $E_j(z)$. We first introduce a structural equation model representation for swing equation dynamics, which will be referred as Linear Dynamic Graph(LDG) for Swing Equations.

DEFINITION 1. *(LDG of Swing Equations) The Linear Dynamic Graph(LDG) of swing equations is defined as $\hat{\mathcal{G}} = (H(z), E(z))$, where, $H(z)$ is the $N \times N$ matrix of transfer functions such that $[H(z)](j,i) = H_{ji}(z), i \neq j$, $[H(z)](j,j) = 0$ and $E(z) = (E_1(z), E_2(z), \cdots, E_N(z))^T$. $\hat{\mathcal{G}}$ corresponds to a directed graph with vertex set $\mathcal{V}$ and edge set $\hat{\mathcal{E}}$, where, $\hat{\mathcal{E}} := \{(i \to j) : [H(z)](j,i) \neq 0\}$.*

Thus the resulting LDG $\hat{\mathcal{G}} = (H(z), E(z))$ has two directed edges for each undirected edge in the grid graph. An illustrative example for the LDG of swing equations $\hat{\mathcal{G}}$ for a given grid graph $\mathcal{G}$ is shown in Fig. 2. We define $\mathcal{N}_j := \{i \in \mathcal{V} : (i, j) \in \mathcal{E}\}$ and $\mathcal{N}_{j,2} := \{i \in \mathcal{V} : (j, k), (i, k) \in \mathcal{E}, \text{for some } k \in \mathcal{V}\}$ as the set of neighbors and two-hop neighbors of $j$ in $\mathcal{G}$ respectively (a node $i$ can be a neighbor as well as a two-hop neighbor of some node $j$).

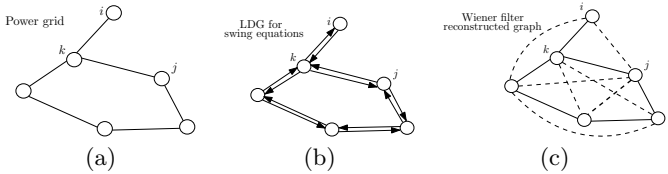

**Figure 2:** (a) undirected Power Grid graph $\mathcal{G}$ (b) bidirected LDG $\hat{\mathcal{G}}$ for swing equations on $\mathcal{G}$ (c) reconstructed graph $\mathcal{G}^w$ obtained using non zero multivariate non causal Wiener filtering.

The 'Wiener filter' or minimum variance estimate of the time series $\theta_j(n)$ by an element $\hat{\theta}_j(n) \in \mathcal{X}_{\bar{j}}$ is given by,

$$\hat{\Theta}_j(z) = (W_{j1}(z) \cdots W_{j(j-1)}(z) W_{j(j+1)}(z) \cdots W_{jN}(z))\Theta_{\bar{j}}(z)$$
$$= \Phi_{\Theta_j \Theta_{\bar{j}}}(z)\Phi_{\Theta_{\bar{j}}}(z)^{-1}\Theta_{\bar{j}}(z), \quad (6)$$

where, $\Theta_{\bar{j}}(z) = (\Theta_1(z) \cdots \Theta_{j-1}(z) \Theta_{j+1}(z) \cdots \Theta_N(z))^T$. The steps involved in computation of an approximate Wiener filter given measurements of nodal dynamics is listed in Appendix A. Moreover, it is shown in [22, 28] that $W_{ji}(z) \neq 0$ implies that $i \in \mathcal{N}_j \cup \mathcal{N}_{j,2}$ (see Theorem B.1 for details). Thus if an undirected graph $\mathcal{G}^w = (\mathcal{V}, \mathcal{E}^w)$ with edge set $\mathcal{E}^w$ is constructed from all non-zero entries of the multivariate non causal Wiener filters $W_{ji}(z)$, then the set of edges will include all edges in the original grid graph $\mathcal{G}$ as well as edges between all two-hop neighbors in $\mathcal{G}$. An example is illustrated in Fig. 2 (c). This leads us to the ***first step*** in our topology learning algorithm where we generate a graph from multivariate non causal Wiener filtering based on the swing dynamics. To determine the grid topology, we need to distinguish between 'true' edges (with neighbors) and 'spurious' edges (with strict two-hop neighbors that are not neighbors). For radial networks, it is possible to distinguish between true and spurious edges due to a local topological separability rule in both static models (see [5]) and [28]). However for loopy networks such topological separability results do not hold in general. In the next section we present a novel pruning algorithm to eliminate the spurious edges obtained by multivariate Wiener filtering.

## 4. PRUNING STEP

The next theorem presents a crucial new result using the phase response of the non causal Wiener filter for spurious edges corresponding to strict two-hop neighbors, which enables us to distinguish them from true edges.

THEOREM 4.1. *In the grid $\mathcal{G} = (\mathcal{V}, \mathcal{E})$, consider a well posed and topologically detectable LDG of swing equations, $\hat{\mathcal{G}} = (H(z), E(z))$. If $i$ and $j$ are strict two-hop neighbors in $\mathcal{G}$ such that $\mathcal{N}_i \cap \mathcal{N}_j \neq \phi$ and $i \notin \mathcal{N}_j, j \notin \mathcal{N}_i$, then $\angle(W_{ji}(e^{\hat{\jmath}\omega})) = -\pi$ for all $\omega \in [-\pi, \pi]$.*

**Proof steps:** Since, $i$ and $j$ are two-hop neighbors, there exist $k \in \mathcal{N}_i \cap \mathcal{N}_j$ such that $b_{kj} > 0$ and $b_{ki} > 0$, implying, $j \to k \leftarrow i$ in the underlying LDG. The proof follows by showing that observing the state at node $k$ induces negative correlation between $i$ and $j$. Aside for pathological parameter cases, the converse of Theorem 4.1 holds making the statement necessary and sufficient for detection of strict two-hop neighbors.

### 4.1 Learning Algorithm

We now present Algorithm 1 that estimates the topology of any general grid $\mathcal{G}$ based on time-series of nodal voltage measurements pertaining to the swing equations. As explained in the preceding sections, it is a two-part algorithm. The first part (Steps 1 - 9) determines the multivariate Wiener filter $W_{ji}(z)$ to estimate the true topology with spurious links between two hop neighbors. In the next part (Steps 10 - 15), we consider a finite set of frequency points $\Omega$ in the interval $[-\pi, \pi)$ and evaluate the phase angle of the Wiener filters for edges in $\mathcal{E}^w$. If the phase angle is within a pre-defined threshold $\tau$ of $-\pi$, the algorithm designates them as spurious edges (see Theorem 4.1) and prunes them from $\mathcal{E}^w$ to produce edge set $\bar{\mathcal{E}}$ of the estimated true topology.

## 5. RESULTS

**Algorithm 1** Topology Learning using Wiener Filtering

**Input:** voltage phase samples $\theta_i$ for nodes $i \in \{1, 2, ...N\}$ in grid $\mathcal{G}$, thresholds $\rho, \tau$, frequency points $\Omega$
**Output:** Estimate of Operational Edges $\bar{\bar{\mathcal{E}}}$

1: **for all** $j \in \{1, 2, ...N\}$ **do**
2:     Compute Wiener filter $W_j(z) = [W_{j1}(z) \cdots W_{jN}(z)]$
3: **end for**
4: Edge set $\mathcal{E}^w \leftarrow \{\}$
5: **for all** $i, j \in \{1, 2, ...N\}, i \neq j$ **do**
6:     **if** $\|W_{ji}(z)\| > \rho$ **then**
7:         $\mathcal{E}^w \leftarrow \mathcal{E}^w \cup \{(i,j)\}$
8:     **end if**
9: **end for**
10: Edge set $\bar{\bar{\mathcal{E}}} \leftarrow \mathcal{E}^w$
11: **for all** $i, j \in \{1, 2, ...N\}, i \neq j$ **do**
12:     **if** $\pi - \tau \leq |\angle(W_{ji}(e^{\hat{\jmath}\omega}))| \leq \pi + \tau, \forall \omega \in \Omega$ **then**
13:         $\bar{\bar{\mathcal{E}}} \leftarrow \bar{\bar{\mathcal{E}}} - \{(i,j)\}$
14:     **end if**
15: **end for**

In this section, we demonstrate the effectiveness of Algorithm 1 presented in the previous section on the IEEE 39 bus test system [2, 11] shown in Fig. 3 (a) with linear dynamics as described by Eq. (3). For our simulations, we model the nodal ambient wide-sense stationary disturbance by white Gaussian noise with spectral density given by $\Phi_{p_j}(z) = 10$ $dB$ and use a small inertia and damping of $0.01$ for nodes without inertia and damping to generate time series data for evaluation of the proposed algorithm. The output at each node is sampled at $0.01s$. To compute the multivariate Wiener filter we use the FIR (Finite Impulse Response) approximation of order $F = 20$ and obtain the filter coefficients by solving linear equations as detailed in Appendix A. For a specific example, consider the the nodes at distance one-hop (colored green) and at distance two-hop (colored red) of node 25 in the IEEE 39 bus system as shown in Fig. 3 (b). The absolute values of the phase response of the multivariate Wiener filters for node 25 and the nodes in its two-hop neighborhood are shown in Fig. 4. It is seen that the phase response of the Wiener filters corresponding to the nodes two hops away are close to $\pi$ rad, while that of the neighbor nodes start from 0 rad. Thus our pruning steps are capable of distinguishing between the two edge types.

For an overall study of the effect of sample size on performance of Algorithm 1, we plot the relative errors in topology estimation for either grid cases in Fig. 5. The threshold $\rho$ was chosen as $10^{-3}$ and $\tau = 0.2\pi$ for the IEEE 39 bus system. It is seen that the relative error, which is defined as the ratio of the sum of false positive and false negative edges to the total number of true edges, decays as the number of samples is increased.

## 6. CONCLUSIONS

In this article, a multivariate Wiener filtering based topology learning approach for the power grid is presented that uses nodal phase angle measurements pertaining to the swing dynamics as input. The Wiener reconstruction of the topology returns spurious links (false positives) between two-hop neighbors in the grid topology along with the true edges.

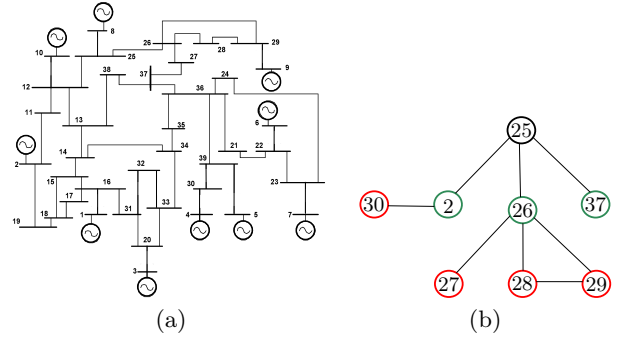

**Figure 3:** (a) IEEE 39 bus system with generators at 10 buses [2], (b) The neighbors (green nodes) and strict two-hop neighbors (red nodes) of node 25 in the IEEE 39 bus system.

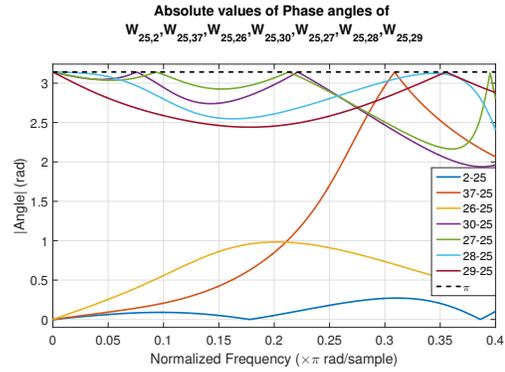

**Figure 4:** Absolute values of the phase response of the Wiener filters between node 25 and its two-hop neighborhood in the IEEE 39 bus system. The phase response begins from 0 rad for all three neighbors and from $\pi$ for all strict 2 hop neighbors of node 25. The length of the time series used is $6.5 \times 10^6$ samples for each node.

The main contribution of our work is in designing a pruning step based on phase response of the Wiener filters to eliminate all spurious links arising out of two-hop neighbor relationship amongst nodes to provably recover the actual topology exactly. Simulation results on IEEE test cases demonstrate the performance of the proposed framework in learning loopy grid topologies.

## 7. ACKNOWLEDGMENTS

The authors S. Talukdar, B. Lundstrom and M. V. Salapaka acknowledge the support of ARPA-E for the project titled 'A Robust Distributed Framework for Flexible Power Grids' via grant no. DE-AR000071. D. Deka acknowledges the support of funding from the U.S. Department of Energy's Office of Electricity as part of the DOE Grid Modernization Initiative.

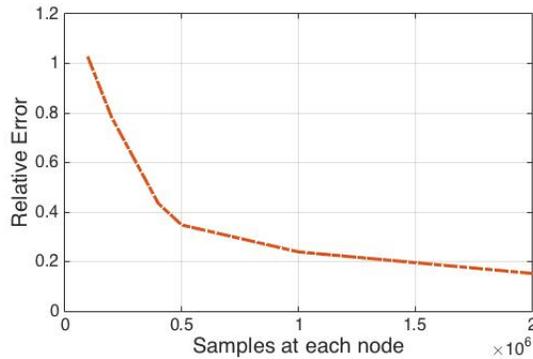

Figure 5: Error proportion with samples per node for IEEE 39 bus system.

# APPENDIX

## A. WIENER FILTER COMPUTATION

As described in Section 3, the Wiener filter determines the optimal projection of a signal $x_j(n)$ in the space $\mathcal{X}_{\bar{j}} = \text{span}\{\theta_i(n+p) : p \in \mathbb{Z}\}_{i \neq j}$. Here we compute the Wiener filter by approximating it with a finite impulse response (FIR) filter, also known as FIR Wiener filter. Let the order of the FIR Wiener filter be $F$. Here the optimal estimate $\hat{\theta}_j(n)$ is written as,

$$\hat{\theta}_j(n) = \sum_{k \in \mathcal{V}, k \neq j} \sum_{p=-F}^{F} h_{k,l} \theta_k(n+l) \quad (7)$$

The Wiener filtering orthogonality condition is stated below and is used to determine the constants $h_{k,l}$ in $\hat{\theta}_j(n)$.

$$\mathbb{E}[\hat{\theta}_j(n)\theta_i(n+l)] = \mathbb{E}[\theta_j(n)\theta_i(n+l)], \quad (8)$$
$$\text{for all } i \in \mathcal{V}, i \neq j, l \in \{-F, -F+1, \cdots, F-1, F\}. \quad (9)$$

Thus, using (7) in (8),

$$[R_{\theta_1\theta_i}(-F-l) \cdots R_{\theta_1\theta_i}(F-l) \cdots R_{\theta_N\theta_i}(-F-l) \cdots R_{\theta_N\theta_i}(F-l)]h$$
$$= R_{\theta_j\theta_i}(-l), i \in \mathcal{V}, i \neq j, l \in \{-F, \cdots, F\}, \quad (10)$$
$$h := [h_1^T \ h_2^T \cdots h_{j-1}^T \ h_{j+1}^T \cdots h_N^T]^T, \text{and}$$
$$h_i^T := [h_{i,-F} \cdots h_{i,-1} \ h_{i,0} \ h_{i,1} \cdots h_{i,F}].$$

The set of equations in (10) describe $(2F+1)(N-1)$ linear equations in $(2F+1)(N-1)$ unknowns in the vector $h$. Thus, in combined form the equations in (10) can be written as,

$$Rh = S,$$

Thus, $h = R^{-1}S$ is used to compute the coefficients of the Wiener filters. Note that the matrix $R$ and the vector $S$ can be computed using empirical correlations from the measured time series data. More details on numerical aspects of Wiener filtering can be found in [14, 16].

## B. WIENER FILTERING AND TOPOLOGY LEARNING

THEOREM B.1. [22, 28] *Consider a LDG of swing equations, $\hat{\mathcal{G}} = (H(z), E(z))$. Let the output of the LDG at the $n^{th}$ time instant be given by $\theta(n) = (\theta_1(n) \ ... \ \theta_N(n))^T$. Define the space $\mathcal{X}_{\bar{j}} = \text{span}\{\theta_i(p) : p \in \mathbb{Z}\}_{i \neq j}$ and in z-domain $\mathcal{X}_{\bar{j}}(z) = \text{span}\{z^p\Theta_i(z) : p \in \mathbb{Z}\}_{i \neq j}$. The 'Wiener filter' of the time series $\theta_j(n)$ by an element $\hat{\theta}_j(n) \in \mathcal{X}_{\bar{j}}$ is given by*

$$\hat{\Theta}_j(z) = \Phi_{\Theta_j \Theta_{\bar{j}}}(z) \Phi_{\Theta_{\bar{j}}}(z)^{-1} \Theta_{\bar{j}}(z), \quad (11)$$

*where, $\Theta_{\bar{j}}(z) = (\Theta_1(z) \ \cdots \ \Theta_{j-1}(z) \ \Theta_{j+1}(z) \ \cdots \ \Theta_N(z))^T$. Also, $W_{ji}(z) \neq 0$ implies that $i \in \mathcal{N}_j \cup \mathcal{N}_{j,2}$.*